
%

\input amstex
\loadbold
\documentstyle{amsppt}
\NoBlackBoxes

\pagewidth{32pc}
\pageheight{44pc}
\magnification=\magstep1

\def\QQ{\Bbb{Q}}
\def\RR{\Bbb{R}}
\def\CC{\Bbb{C}}
\def\PP{\Bbb{P}}
\def\OO{\Cal{O}}
\def\Spec{\operatorname{Spec}}
\def\zeros{\operatorname{div}}
\def\rank{\operatorname{rk}}
\def\codim{\operatorname{codim}}
\def\rest#1#2{\left.{#1}\right\vert_{{#2}}}
\def\submet#1#2{\rest{#1}{#2}}
\def\achern#1#2{\widehat{c}_{#1}(#2)}
\def\achow{\widehat{\operatorname{CH}}}
\def\chow{\operatorname{CH}}
\def\adeg{\widehat{\operatorname{deg}}}
\def\poscone{\widehat{C}_{+\kern-.09em+}}
\def\wposcone{\widehat{C}_{+}}
\def\adisc{\widehat{\delta}}
\def\adiff{\widehat{d}}
\def\Gal{\operatorname{Gal}}

\topmatter
\title
Bogomolov unstability on arithmetic surfaces
\endtitle
\rightheadtext{}
\author Atsushi Moriwaki \endauthor
\leftheadtext{}
\address
Department of Mathematics, Faculty of Science,
Kyoto University, Kyoto, 606-01, Japan
\endaddress
\curraddr
Department of Mathematics, University of California,
Los Angeles, 405 Hilgard Avenue, Los Angeles, California 90024, USA
\endcurraddr
\email moriwaki\@math.ucla.edu \endemail
\date March, 1994 \enddate
\abstract
In this paper, we will consider an arithmetic analogue
of Bogomolov unstability theorem, i.e.
if $(E, h)$ is a torsion free Hermitian sheaf on an arithmetic
surface $X$ and
$\adeg\left((\rank E - 1) \achern{1}{E, h}^2 - (2 \rank E)
\achern{2}{E, h}\right) > 0$,
then there is a non-zero saturated subsheaf $F$ of $E$ such that
$\achern{1}{F, \submet{h}{F}}/\!\rank F - \achern{1}{E, h}/\!\rank E$
lies in the positive cone of $X$.
\endabstract
\endtopmatter

\document

\head
0. Introduction
\endhead

In \cite{Bo}, Bogomolov proved unstability theorem, namely,
if a vector bundle $E$ on a complex projective surface $S$
satisfies an inequality
$(\rank E - 1) c_1(E)^2 - (2 \rank E) c_2(E) > 0$, then
there is a saturated subsheaf $F$ of $E$ such that
$c_1(F)/\!\rank F - c_1(E)/\!\rank E$ belongs to the positive cone
of $S$.
In this paper, we would like to consider an arithmetic analogue
of Bogomolov unstability theorem.

Let $f : X \to \Spec(O_K)$ be a regular arithmetic surface
over the ring of integers of a number field $K$ with $f_*\OO_X = O_K$,
and $\deg_K  : \achow^1(X)_{\RR} \to \RR$ the natural homomorphism
given by
$$
\deg_K : \achow^1(X)_{\RR} \overset z\to\longrightarrow
\chow^1(X)_{\RR} \overset\otimes K\to\longrightarrow
\chow^1(X_K)_{\RR} \overset \deg\to\longrightarrow \RR.
$$
The positive cone $\poscone(X)$ of $X$
is defined by the set of all elements $x \in \achow^1(X)_{\RR}$
with $\adeg(x^2) > 0$ and $\deg_K(x) > 0$.
A torsion free Hermitian sheaf $(E, h)$ on $X$
is said to be {\it arithmetically unstable} if
there is a non-zero saturated subsheaf $F$ of $E$ with
$$
\frac{\achern{1}{F, \submet{h}{F}}}{\rank F} - \frac{\achern{1}{E, h}}{\rank E}
\in \poscone (X),
$$
where $\submet{h}{F}$ is the Hermitian metric of $F$ given by
the restriction of $h$ to $F$.
The main theorem of this paper is the following.

\proclaim{Theorem A}
If $(E, h)$ is a torsion free Hermitian sheaf on $X$ and
$$
\adeg\left((\rank E - 1) \achern{1}{E, h}^2 - (2 \rank E)
\achern{2}{E, h}\right) > 0,
$$
then $(E, h)$ is arithmetically unstable.
\endproclaim

Let $A$ be an element of $\achow^1(X)_{\RR}$ with
$\adeg(A \cdot x) > 0$ for all $x \in \poscone(X)$, i.e.
according to terminology in \S1.1,
$A$ is an element of the weak positive cone $\wposcone(X)$.
A torsion free Hermitian sheaf $(E, h)$ on $X$ is said to be
{\it arithmetically $\mu$-semistable with respect to $A$}
if, for all non-zero saturated subsheaves $F$ of $E$,
$$
\frac{\adeg(\achern{1}{F, \submet{h}{F}} \cdot A)}{\rank F} \leq
\frac{\adeg(\achern{1}{E, h} \cdot A)}{\rank E}.
$$
Then, we have the following corollary of Theorem~A.

\proclaim{Corollary~B}
If $(E, h)$ is arithmetically $\mu$-semistable with respect to $A$,
then
$$
\adeg\left((2 \rank E) \achern{2}{E, h} -
(\rank E - 1)\achern{1}{E, h}^2
\right) \geq 0.
$$
\endproclaim

If $A = \left(0, \sum_{\sigma \in K(\CC)} 1/[K : \QQ]\right)$, then
arithmetic $\mu$-semistability of $(E, h)$ with respect to $A$
is nothing more than
$\mu$-semistability of $E_{\overline{\QQ}}$.
So this corollary gives a generalization of \cite{Mo1}, \cite{Mo2}
and \cite{So}.

In \S1, we will prepare several basic facts of the positive cone
and Hermitian vector spaces.
In \S2, we will consider finiteness of saturated subsheaves
in a Hermitian vector bundle, which will be crucial
for the proof of Theorem~A.
\S3 is devoted to the proof of Theorem~A and Corollary~B.

\head
1. Preliminaries
\endhead

\subhead
1.1. Positive cone of arithmetic Chow group
\endsubhead
Here, we consider basic properties
of the positive cone of the arithmetic Chow group of
codimension $1$.

Let $K$ be a number field and $O_K$ the ring of integers of $K$.
Let $f : X \to \Spec(O_K)$ be a regular arithmetic surface
with $f_*\OO_X = O_K$,
and $\deg_K  : \achow^1(X)_{\RR} \to \RR$ the natural homomorphism
defined by
$$
\deg_K : \achow^1(X)_{\RR} \overset z\to\longrightarrow
\chow^1(X)_{\RR} \overset\otimes K\to\longrightarrow
\chow^1(X_K)_{\RR}
\overset \deg\to\longrightarrow \RR.
$$
We set
$$
\align
\poscone(X) & = \left\{ x \in \achow^1(X)_{\RR} \mid
\hbox{$\adeg(x^2) > 0$ and $\deg_K(x) > 0$}
\right\} \quad\hbox{and} \\
\wposcone(X)  & = \left\{ x \in \achow^1(X)_{\RR} \mid
\hbox{$\adeg(x \cdot y) > 0$ for all $y \in \poscone(X)$} \right\}.
\endalign
$$
$\poscone(X)$ (resp. $\wposcone(X)$)
is called {\it the positive cone of $X$}
(resp. {\it the weak positive cone of $X$}).
First of all, we have the following lemma.

\proclaim{Lemma 1.1.1}
\rom{(1.1.1.1)} If $h \in \poscone(X)$, $x \in \achow^1(X)_{\RR}$
and $\adeg(x \cdot h) = 0$, then $\adeg(x^2) \leq 0$.

\rom{(1.1.1.2)} If $x \in \achow^1(X)_{\RR}$,
$\adeg(x^2) \geq 0$ and $\deg_K(x) > 0$,
then $x \in \wposcone(X)$.

\rom{(1.1.1.3)} If $h \in \poscone(X)$, $x \in \achow^1(X)_{\RR}$,
$\adeg(x^2) \geq 0$ and $\adeg(x \cdot h) > 0$,
then $x \in \wposcone(X)$.
\endproclaim

\demo{Proof}
(1.1.1.1) Let $t$ be a real number with
$\deg_K(x - th) = 0$. Then,
by Hodge index theorem (cf. \cite{Fa}, \cite{Hr} or \cite{Mo3}),
$\adeg((x - t h)^2) \leq 0$. Thus, $\adeg(x^2) + t^2 \adeg(h^2) \leq 0$.
Therefore, $\adeg(x^2) \leq 0$.

(1.1.1.2)
Let $y \in \poscone(X)$ and $t$ a real number with
$\deg_K(y - t x) = 0$. Then, $t > 0$.
By Hodge index theorem, $\adeg((y - t x)^2) \leq 0$.
Thus, we have
$$
   \adeg(y^2) + t^2 \adeg(x^2) \leq 2t \adeg(x \cdot y).
$$
Therefore, $\adeg(x \cdot y) > 0$.
Hence, $x \in \wposcone(X)$.

(1.1.1.3)
Let $y \in \poscone(X)$ and $t$ a real number with
$\adeg(y - tx \cdot h) = 0$. Then, $t > 0$ by (1.1.1.2).
(1.1.1.1) implies $\adeg((y - t x)^2) \leq 0$.
Thus, by the same way as above,
we have $\adeg(x \cdot y) > 0$, which implies
$x \in \wposcone(X)$.
\qed
\enddemo

$\poscone(X)$ and $\wposcone(X)$ have the following properties.

\proclaim{Proposition 1.1.2}
\rom{(1.1.2.1)} $\poscone(X) \subset \wposcone(X)$.

\rom{(1.1.2.2)} If $x, y \in \poscone(X)$ and $t > 0$,
then $x+y, tx \in \poscone(X)$.

\rom{(1.1.2.3)} If $x, y \in \wposcone(X)$ and $t > 0$,
then $x+y, tx \in \wposcone(X)$.

\rom{(1.1.2.4)} $\poscone(X) =
\left\{ x \in \achow^1(X)_{\RR} \mid
\hbox{$\adeg(x \cdot y) > 0$ for all
$y \in \wposcone(X)$} \right\}$.
\endproclaim

\demo{Proof}
(1.1.2.1) and (1.1.2.2) are straightforward from (1.1.1.2).
(1.1.2.3) is obvious.

(1.1.2.4) Clearly, we have
$$
\poscone(X) \subseteq
\left\{ x \in \achow^1(X)_{\RR} \mid
\hbox{$\adeg(x \cdot y) > 0$ for all
$y \in \wposcone(X)$} \right\}.
$$
We assume that $\adeg(x \cdot y) > 0$ for all
$y \in \wposcone (X)$.
If we set $B = \left(0, \sum_{\sigma \in K(\CC)} 1/[K : \QQ]\right)$,
then $B \in \wposcone(X)$. Thus, $\deg_K(x) = \adeg(x \cdot B) > 0$.
Hence, it is sufficient to show $\adeg(x^2) > 0$.
Here, we fix $h \in \poscone(X)$ with $\adeg(h^2) = 1$.
We set $t = \adeg(x \cdot h) > 0$.
Since $\adeg(x - th \cdot h) = 0$, by (1.1.1.1),
$\adeg((x - th)^2) \leq 0$. If $\adeg((x - th)^2) = 0$,
then
$$
\adeg(x^2) = \adeg ((th + (x - th))^2) = t^2 > 0.
$$
Thus, we may assume that $\adeg((x - th)^2) < 0$
If we set
$s = \left(-\adeg((x -  th)^2)\right)^{1/2}$ and
$l = (x - th)/s$,
then $x = th + sl$, $\adeg(l^2) = -1$ and
$\adeg(h \cdot l) = 0$.
Let us consider $y = h + l$.
Since $\adeg(y^2) = 0$ and $\adeg(y \cdot h) = 1$,
by (1.1.1.3), $y \in \wposcone (X)$.
Thus, $\adeg(x \cdot y) = t - s > 0$.
Therefore
$$
\adeg(x^2) = t^2 - s^2 = (t + s)(t - s) > 0.
$$
Hence, $x \in \poscone (X)$.
\qed
\enddemo

Finally we consider the following proposition.

\proclaim{Proposition 1.1.3}
For $z \in \achow^1(X)_{\RR}$, we set
$$
W(z) = \left\{ u \in \wposcone(X) \mid
\adeg(z \cdot u) > 0 \right\}.
\tag 1.1.3.1
$$
If $x \not\in \poscone(X)$,
$y \in \poscone(X)$ and $W(x) \not= \emptyset$,
then $W(x) \subsetneqq  W(x + y)$.
\endproclaim

\demo{Proof}
Clearly, $W(x) \subseteq W(x+y)$.
By virtue of (1.1.2.4),
$W(x) \subsetneqq \wposcone(X)$.
Let $u_1 \in W(x)$, $u_2 \in \wposcone(X) \setminus W(x)$
and $t = -\adeg(x \cdot u_2)/\adeg(x \cdot u_1)$. Then,
$t \geq 0$ and $\adeg(x \cdot u_2 + t u_1) = 0$.
Hence, $u_2 + t u_1 \not\in W(x)$.
On the other hand, by (1.1.2.3), $u_2 + t u_1 \in \wposcone(X)$.
Moreover, $\adeg(y \cdot u_2 + t u_1) > 0$.
Thus, $\adeg(x + y \cdot u_2 + t u_1) > 0$.
Therefore $u_2 + t u_1 \in W(x+y)$.
\qed
\enddemo

\subhead
1.2. Hermitian vector space
\endsubhead
Let $V$ be a $\CC$-vector space, $h_V$ a Hermitian metric on $V$ and
$W$ a subvector space of $V$. Considering the restriction
of $h$ to $W$, the metric $h_V$ induces a metric $h_W$ of
$W$, which is called {\it the submetric of $W$ induced by $h_V$}.
Let $W^{\perp}$ be the orthogonal complement of $W$.
Then the natural homomorphism $W^{\perp} \longrightarrow V/W$ is isomorphic.
Thus we have a metric $h_{V/W}$ of $V/W$ given by $\rest{h}{W^{\perp}}$.
This metric is called
{\it the quotient metric of $V/W$ induced by $h_V$}.

\proclaim{Proposition 1.2.1}
Let $(V, h_V)$ be a Hermitian vector space over $\CC$ and
$U, W$ subspaces of $V$ with $U \subset W$.
Let $h_W = \rest{(h_V)}{W}$ and $h_{V/U}$ the quotient
metric of $V/U$ induced by $h_V$.
We consider two Hermitian metrics of $W/U$.
Let $h_s = \rest{(h_{V/U})}{W/U}$ and $h_q$
the quotient metric induced by $h_W$. Then, we have $h_s = h_q$.
\endproclaim

\demo{Proof}
Let $V = U \oplus U^{\perp}$ be the orthogonal decomposition of $V$
and $f : U^{\perp} \to V/U$ the natural isomorphism.
Then, for $x, y \in W/U$, $h_s(x, y) = h_V(f^{-1}(x), f^{-1}(y))$.

Let $W = U \oplus U_W^{\perp}$ be the orthogonal decomposition of $W$,
i.e. $U_W^{\perp}$ is the orthogonal complement of $U$ in $W$,
and $g : U_W^{\perp} \to W/U$ the natural isomorphism.
Then, $h_q(x, y) = h_V(g^{-1}(x), g^{-1}(y))$ for $x, y \in W/U$.

On the other hand, $U_W^{\perp} \subset U^{\perp}$ and the
following diagram is commutative.
$$
\CD
U_W^{\perp} & {\ }\hookrightarrow{\ } &   U^{\perp} \\
@V{g}VV       @VV{f}V   \\
W/U         & {\ }\hookrightarrow{\ } &   V/U
\endCD
$$
Thus, we have $h_s = h_q$.
\qed
\enddemo

\head
2. Finiteness of saturated subsheaves
\endhead

In this section, we will consider finiteness of saturated subsheaves
in a Hermitian vector bundle, which will be crucial
for the proof of Theorem~A.

\proclaim{Theorem 2.1}
Let $K$ be a number field and $O_K$ the ring of integers of $K$.
Let $f : X \to \Spec(O_K)$ be a regular arithmetic surface
with $f_*\OO_X = O_K$,
$(E, h)$ a Hermitian vector bundle on $X$ and
$(H, k)$ a Hermitian line bundle on $X$. If $H_K$ is ample,
then, for constants $C_1$ and $C_2$, the set
$$
\Cal{F} = \left\{ L \ \left| \
\aligned
& \hbox{$L$ is a rank-1 saturated subsheaf of $E$ with} \\
& \hbox{$\adeg(\achern{1}{L, \submet{h}{L}} \cdot \achern{1}{H, k}) \geq C_1$
and $\deg(L_K) \geq C_2$}
\endaligned
\right.\right\}
$$
is finite.
\endproclaim

\demo{Proof}
First of all, we need the following Lemma.

\proclaim{Lemma 2.1.1}
Let $K$ be an infinite field, $C$ a smooth projective curve over $K$, and
$E$ a vector bundle on $C$ of rank $r \geq 2$.
Then, for every real number $M$, there is a rank-1 saturated subsheaf
$L$ of $E$ with $\deg(L) < M$.
\endproclaim

\demo{Proof}
Let $H$ be an ample line bundle on $X$ and
$n$ a positive integer such that
$n \deg(H) > -M$ and $E \otimes H^n$ is generated by global sections.
Let us consider the following closed set $\Sigma$ in
$C \times \PP(H^0(C, E \otimes H^n))$.
$$
\Sigma = \left\{ (x, s) \in C \times \PP(H^0(C, E \otimes H^n))
\mid s(x) = 0 \right\}.
$$
Let $p : \Sigma \to C$ be the natural projection.
Since $E \otimes H^n$ is generated by global sections,
$\codim(p^{-1}(x), \PP(H^0(C, E \otimes H^n))) = r$.
Therefore, $\codim(\Sigma, C \times \PP(H^0(C, E \otimes H^n)))= r$,
which means $\dim \Sigma < \dim \PP(H^0(C, E \otimes H^n))$.
Thus, the natural projection $q : \Sigma \to \PP(H^0(C, E \otimes H^n))$
is not surjective.
Hence, since $\#(K)$ is infinite,
$$
  (p(\Sigma))(K) \subsetneqq \PP(H^0(C, E \otimes H^n))(K).
$$
Thus, there is a section $s \in H^0(C, E \otimes H^n)$ with
$\zeros(s) = \emptyset$, which induces an injective homomorphism
$H^{-n} \to E$. Since $\zeros(s) = \emptyset$,
the image of $H^{-n} \to E$ is saturated.
\qed
\enddemo

Let us start the proof of Theorem~2.1.
Clearly we may assume $r = \rank E \geq 2$.
By Lemma~2.1.1, there is a filtration :
$\{ 0 \} = F_0 \subset F_1 \subset \cdots \subset F_{r-1}
\subset F_r = E$ such that

(2.1.2) $F_i/F_{i-1}$ is a rank-1 torsion free sheaf
for every $1 \leq i \leq r$.

(2.1.3) $\deg((F_i/F_{i-1})_K) < C_2$ for
$ 1 \leq i \leq r-1$.

\noindent
Here we claim

\proclaim{Claim 2.1.4} If $L$ is a line bundle
on $X$ with $\deg(L_K) \geq C_2$ and
$\varphi : L \to F_{r-1}$ is a homomorphism,
then $\varphi = 0$.
\endproclaim

We assume that $\varphi \not= 0$.
Choose $i$ in such a way that $\varphi(L) \subseteq F_i$ and
$\varphi(L) \not\subseteq F_{i-1}$. Then,
we have an injective homomorphism
$L \to F_i/F_{i-1}$. Since $i \leq r-1$,
$\deg(L_K) \geq C_2 > \deg((F_i/F_{i-1})_K)$.
This is a contradiction.
\qed

\bigskip
Let $Q$ be the double dual of $E/F_{r-1}$ and $h_Q$
the quotient metric of $Q$ induced by $h$ via $E \to E/F_{r-1}$, i.e.
$h_Q = \submet{h}{F^{\perp}}$.
Pick up $L \in \Cal{F}$.
By Claim~2.1.4, we have the natural injection
$L \to Q$.  So there is an effective divisor $D_L$ on $X$
such that $L \otimes \OO_X(D_L) \simeq Q$.
Let $D^h_L$ (resp. $D^v_L$) be the horizontal part of $D_L$
(resp. the vertical part of $D_L$). Then, we have

\proclaim{Claim 2.1.5}
If $D^h_L = D^h_{L'}$ for $L, L' \in \Cal{F}$,
then $L = L'$.
\endproclaim

Let us consider $M = L \otimes \OO_X(-D^v_{L'})$ and
$M' = L' \otimes \OO_X(-D^v_{L})$. Since $D^h_L = D^h_{L'}$,
$M$ and $M'$ has the same image in $Q$ via $E \to Q$.
Therefore, we have
$M + F_{r-1} = M' + F_{r-1}$. Moreover,
since $M \to Q$ and $M' \to Q$ are injective,
$M \cap F_{r-1} = M' \cap F_{r-1} = \{ 0 \}$.
Hence, we have a homomorphism
$M \to M' \oplus F_{r-1} \to F_{r-1}$,
which must be zero by Claim~2.1.4.
Thus $M \subseteq M'$. By the same way, $M' \subseteq M$.
Hence, $M = M'$. So $L = L'$ because
$L$ is the saturation of $M$ in $E$ and
$L'$ is the saturation of $M'$ in $E$.
\qed

\bigskip
We set $C_3 = \adeg(\achern{1}{Q, h_Q} \cdot \achern{1}{H, k}) - C_1$ and
$C_4 = \deg(Q_K) - C_2$.
Let $D^h_L = \sum_i a_i \Gamma_i$ be the irreducible
decomposition of $D^h_L$. Then, we have

\proclaim{Claim~2.1.6}
$\sum_i a_i \adeg(\achern{1}{\rest{(H, k)}{\Gamma_i}}) \leq C_3$ and
$\sum_i a_i [K(\Gamma_i) : K] \leq C_4$.
\endproclaim

Since $L \otimes \OO_X(D_L) \simeq Q$, we have
$\deg(L_K) + \deg((D_L)_K) = \deg(Q_K)$. Thus we get
$\sum_i a_i [K(\Gamma_i) : K] \leq C_4$.

We choose a Hermitian metric $h_{D_L}$ of $\OO_X(D_L)$
in such a way that $(L, \submet{h}{L}) \otimes (\OO_X(D_L), h_{D_L})$
is isometric to $(Q, h_Q)$. Let $1$ be the canonical section
of $H^0(X, \OO_X(D_L))$ with $\zeros(1) = D_L$, and
$D^v_L = \sum_j b_j l_j$ the irreducible decomposition of the
vertical part of $D_L$. Then, we have
$$
\multline
\adeg(\achern{1}{\OO_X(D_L), h_{D_L}} \cdot \achern{1}{H, k})  =
      \sum_i a_i \adeg(\achern{1}{\rest{(H, k)}{\Gamma_i}}) +
      \sum_j b_j \deg(\rest{H}{l_j}) \\
 - \frac{1}{2} \sum_{\sigma \in K(\CC)}
      \int_{X_{\sigma}} \log(h_{D_L}(1, 1)) c_1(H_{\sigma}, k_{\sigma}).
\endmultline
$$
Since $h_Q$ is a quotient metric of $h$, we can see
$h_{D_L}(1, 1) \leq 1$ for all points of each infinite fiber $X_{\sigma}$.
Therefore, we get
$$
\align
\sum_i a_i \adeg(\achern{1}{\rest{(H, k)}{\Gamma_i}}) & \leq
\adeg(\achern{1}{\OO_X(D_L), h_{D_L}} \cdot \achern{1}{H, k}) \\
& = \adeg(\achern{1}{Q, h_Q} \cdot \achern{1}{H, k})
- \adeg(\achern{1}{L, \submet{h}{L}} \cdot \achern{1}{H, k}) \\
& \leq C_3.
\endalign
$$
Thus, we obtain our claim.
\qed

\bigskip
To complete our theorem, by Claim~2.1.5, it is sufficient to see that
$\{ D^h_L \mid L \in \Cal{F} \}$ is finite.
Since $\sum_i a_i [K(\Gamma_i) : K] \leq C_4$, we have
$a_i \leq C_4$ and $[\QQ(\Gamma_i) : \QQ] \leq C_4 [K : \QQ]$.
Let $D^h_L = D^h_L(+) + D^h_L(-)$ be the decomposition of $D^h_L$
such that
$$
D^h_L(+) = \sum_{\adeg(\achern{1}{\rest{(H, k)}{\Gamma_i}}) > 0}
a_i \Gamma_i
\quad\hbox{and}\quad
D^h_L(-) = \sum_{\adeg(\achern{1}{\rest{(H, k)}{\Gamma_i}}) \leq 0}
a_i \Gamma_i.
$$
By Northcott's theorem (cf. Theorem 2.6 in Chapter 2 of \cite{La}),
the set $\{ D^h_L(-) \mid L \in \Cal{F} \}$ is finite.
Hence, there is a constant $C_5$ depending only on $\Cal{F}$ with
$$
\sum_{\adeg(\achern{1}{\rest{(H, k)}{\Gamma_i}}) > 0}
 a_i \adeg(\achern{1}{\rest{(H, k)}{\Gamma_i}}) \leq C_5.
$$
Thus, $\{ D^h_L(+) \mid L \in \Cal{F} \}$ is finite.
Therefore, $\{ D^h_L \mid L \in \Cal{F} \}$ is finite.
\qed
\enddemo

\proclaim{Corollary 2.2}
Let $K$ be a number field and $O_K$ the ring of integers of $K$.
Let $f : X \to \Spec(O_K)$ be a regular arithmetic surface
with $f_*\OO_X = O_K$,
$(E, h)$ a Hermitian vector bundle on $X$ and
$(H, k)$ a Hermitian line bundle on $X$.
If $H_K$ is ample, then,
for constants $C_1$ and $C_2$, the set
$$
\left\{ \achern{1}{F, \submet{h}{F}} \in \achow^1(X) \ \left| \
\aligned
& \hbox{$F$ is a non-zero saturated subsheaf of $E$ with} \\
& \hbox{$\adeg(\achern{1}{F, \submet{h}{F}} \cdot \achern{1}{H, k}) \geq C_1$
and $\deg(F_K) \geq C_2$}
\endaligned
\right.\right\}
$$
is finite.
\endproclaim

\demo{Proof}
Since $\det F$ is a saturated subsheaf of
${\displaystyle \bigwedge^{\rank F} E}$,
our corollary is an immediate consequence of Theorem~2.1.
\qed
\enddemo

\head
3. Proof of Bogomolov unstability (Theorem A)
\endhead

Before the proof of Theorem~A, we will fix notations.
Let $K$ be a number field and $O_K$ the ring of integers of $K$.
Let $f : X \to \Spec(O_K)$ be a regular arithmetic surface
with $f_*\OO_X = O_K$.
Let $(F, h_F)$ and $(E, h_E)$ be torsion free Hermitian
sheaves on $X$. We set
$$
 \adisc(E, h_E) =
 \adeg \left(
 \frac{\rank E - 1}{2 \rank E} \achern{1}{E, h_E}^2 - \achern{2}{E, h_E}
 \right).
$$
and
$$
\adiff((F, h_F), (E, h_E)) =
\frac{\achern{1}{F, h_F}}{\rank F} - \frac{\achern{1}{E, h_E}}{\rank E}.
$$
Then, we have the following formula.

\proclaim{Lemma 3.1}
Let $0 \to (S, h_S) \to (E, h_E) \to (Q, h_Q) \to 0$
be an exact sequence of torsion free Hermitian sheaves on $X$ such that
$h_S$ and $h_Q$ are the induced metric by $h_E$.
Then, we have
$$
 \adisc(E, h_E) \leq \adisc(S, h_S) + \adisc(Q, h_Q) +
 \frac{(\rank E)(\rank S)}{2 \rank Q}\adeg(\adiff((S, h_S), (E, h_E))^2).
$$
\endproclaim

\demo{Proof}
First of all,
$\achern{1}{E, h_E} = \achern{1}{S, h_S} + \achern{1}{Q, h_Q}$.
Moreover, by Proposition~7.3 of \cite{Mo1},
$$
\adeg\left( \achern{2}{E, h_E} - \achern{2}{(S, h_S) \oplus (Q, h_Q)} \right)
\geq 0.
$$
Thus, by an easy calculation, we have our lemma.
\qed
\enddemo

\bigskip
Let us start the proof of Theorem~A.
Let $E^{\vee\vee}$ be the double dual of $E$.
Then, $\achern{2}{E^{\vee\vee}, h} = \achern{2}{E, h} -
\log\left(\#\left(E^{\vee\vee}/E\right)\right)$.
So we may assume that $E$ is locally free.
First we claim

\proclaim{Claim 3.2}
There is a non-zero saturated subsheaf $F$ of $E$ such that
$$
\frac{\deg(F_K)}{\rank F} - \frac{\deg(E_K)}{\rank E} > 0.
$$
\endproclaim

Since
$\adeg\left((\rank E - 1) \achern{1}{E, h}^2 - (2 \rank E)
\achern{2}{E, h}\right) > 0$,
by the main theorem of \cite{Mo1}, $E_{\overline{\QQ}}$ is not semistable.
Let $F'$ be the maximal destabilizing sheaf of $E_{\overline{\QQ}}$.
For $\tau \in \Gal(\overline{\QQ}/K)$, let us consider $\tau(F')$.
Then, $\tau(F') \subset E_{\overline{\QQ}}$,
$\deg(\tau(F')) = \deg(F')$ and $\rank \tau(F') = \rank F'$,
which means that $\tau(F')$ is also a maximal destabilizing sheaf of
$E_{\overline{\QQ}}$. Thus, by the uniqueness of the
maximal destabilizing sheaf, we have $\tau(F') = F'$.
Therefore, $F'$ is defined over $K$.
Hence, there is a saturated subsheaf $F$ of $E$ with $F_K = F'$.
Thus, we have our claim.
\qed

\bigskip
Let $(H, k)$ be a Hermitian line bundle on $X$ such that
$H_K$ is ample.
Since
$$
\multline
	\adeg(\adiff((F, \submet{h}{F}), (E, h)) \cdot \achern{1}{H, ck})
= \adeg(\adiff((F, \submet{h}{F}), (E, h)) \cdot \achern{1}{H, k}) \\
- \frac{\log(c)[K : \QQ]}{2} \left(
\frac{\deg(F_K)}{\rank F} - \frac{\deg(E_K)}{\rank E} \right)
\endmultline
$$
and $\adeg(\achern{1}{H, ck}^2) = \adeg(\achern{1}{H, k}^2) -
\log(c)[K : \QQ]\deg(H_K)$,
we may assume that
$\adeg(\adiff((F, \submet{h}{F}), (E, h))) \cdot \achern{1}{H, k}) > 0$
and $\achern{1}{H, k} \in \poscone(X)$
if we replace $k$ by $ck$ with sufficiently small positive number $c$.
Here we consider the following set.
$$
\Cal{G} = \left\{  G \ \left| \
\aligned
& \hbox{$G$ is a non-zero saturated subsheaf of $E$ with} \\
& \hbox{$\adeg(\adiff((G, \submet{h}{G}), (E, h)) \cdot \achern{1}{H, k}) > 0$
and} \\
& \hbox{$\deg_K(\adiff((G, \submet{h}{G}), (E, h))) > 0$}
\endaligned
\right.\right\}.
$$
Then, $F \in \Cal{G}$. Moreover,
by Corollary~2.2, $\left\{
\adiff((G, \submet{h}{G}), (E, h)) \in \achow^1(X)_{\QQ} \mid
G \in \Cal{G} \right\}$ is finite.

We will prove Theorem~A by induction on $\rank E$.

\proclaim{Claim 3.3}
If $\rank E = 2$, then
	$\adiff((F, \submet{h}{F}), (E, h)) \in \poscone(X)$.
\endproclaim

Let $h_{E/F}$ be the quotient metric of $E/F$, i.e.
$h_{E/F} = \submet{h}{F^{\perp}}$.
By Lemma~3.1,
$$
 \adisc(E, h) \leq \adisc(F, \submet{h}{F}) + \adisc(E/F, h_{E/F}) +
 \frac{(\rank E)(\rank F)}{2 \rank E/F}\adeg(\adiff((F, \submet{h}{F}), (E,
h))^2).
$$
Since $\rank F = \rank E/F = 1$,
$\adisc(F, \submet{h}{F}) \leq 0$ and $\adisc(E/F, h_{E/F}) \leq 0$.
Therefore,
$$
\adeg(\adiff((F, \submet{h}{F}), (E, h))^2) > 0.
$$
Thus, $\adiff((F, \submet{h}{F}), (E, h)) \in \poscone(X)$.
\qed

\bigskip
{}From now on, we assume $\rank E \geq 3$.
As in (1.1.3.1), for $x \in \achow^1(X)_{\RR}$, we set
$$
W(x) = \left\{ u \in \wposcone(X) \mid
\adeg(x \cdot u) > 0 \right\}.
$$
Here we claim

\proclaim{Claim 3.4}
Under the hypothesis of induction,
if $\adeg(\adiff((G, \submet{h}{G}), (E, h))^2) \leq 0$
for $G \in \Cal{G}$,
then there is $G_1 \in \Cal{G}$ with
$W(\adiff((G, \submet{h}{G}), (E, h))) \subsetneqq
W(\adiff((G_1, \submet{h}{G_1}), (E, h)))$.
\endproclaim

We set $h_{E/G} = \submet{h}{G^{\perp}}$.
First of all, by Lemma~3.1,
$$
 \adisc(E, h) \leq \adisc(G, \submet{h}{G}) + \adisc(E/G, h_{E/G}) +
 \frac{(\rank E)(\rank G)}{2 \rank E/G}\adeg(\adiff((G, \submet{h}{G}), (E,
h))^2).
$$
Since $\adisc(E, h) > 0$ and $\adeg(\adiff((G, \submet{h}{G}), (E, h))^2) \leq
0$,
we have either
$\adisc(G, \submet{h}{G}) > 0$ or $\adisc(E/G, h_{E/G}) > 0$.

If $\adisc(G, \submet{h}{G}) > 0$, then by hypothesis of induction
there is a non-zero saturated subsheaf $G_1$ of $G$ with
$\adiff((G_1, \submet{h}{G_1}), (G, \submet{h}{G})) \in \poscone(X)$. Here
since
$$
\adiff((G_1, \submet{h}{G_1}), (E, h)) =
\adiff((G_1, \submet{h}{G_1}), (G, h_G)) + \adiff((G, \submet{h}{G}), (E, h)),
$$
we have $G_1 \in \Cal{G}$. Moreover, by Proposition~1.1.3,
we get
$$
 W(\adiff((G, \submet{h}{G}), (E, h))) \subsetneqq
 W(\adiff((G_1, \submet{h}{G_1}), (E, h))).
$$

If $\adisc(E/G, h_{E/G}) > 0$, then by hypothesis of induction
there is a non-zero saturated subsheaf $T$ of $E/G$ such that
$\adiff((T, \submet{h_{E/G}}{T}), (E/G, h_{E/G})) \in \poscone(X)$.
Take a saturated subsheaf $G_1$ in $E$ with
$G \subset G_1$ and $G_1/G = T$.
Let $h_{G_1/G}$ be the induced quotient metric by $h_{G_1}$.
By Proposition~1.2.1, we have $\submet{h_{E/G}}{T} = h_{G_1/G}$.
Thus, by an easy calculation, we get
$$
\multline
\adiff((G_1, \submet{h}{G_1}), (E, h)) =
\frac{\rank (G_1/G) }{\rank G_1}
\adiff((T, \submet{h_{E/G}}{T}), (E/G, h_{E/G})) \\
             + \frac{\rank G \rank (E/G_1)}
	                  {\rank G_1 \rank (E/G) } \adiff((G, \submet{h}{G}), (E,
h_E)).
\endmultline
$$
Therefore, $G_1 \in \Cal{G}$, and by Proposition~1.1.3,
$$
W(\adiff((G, \submet{h}{G}), (E, h))) \subsetneqq
W(\adiff((G_1, \submet{h}{G_1}), (E, h))).
$$
Hence we get Claim~3.4.
\qed

\bigskip
Here we assume that
$\adeg(\adiff((G, \submet{h}{G}), (E, h))^2) \leq 0$ for all $G \in \Cal{G}$.
Then, since $F \in \Cal{G}$, by Claim~3.4,
there is an infinite sequence $\{G_0 = F, G_1, G_2, \ldots, G_n, \ldots \}$
in $\Cal{G}$ such that
$$
W(\adiff((G_i, \submet{h}{G_i}), (E, h))) \subsetneqq
W(\adiff((G_j, \submet{h}{G_j}), (E, h)))
$$
for all $i < j$.
In particular, $\adiff((G_i, \submet{h}{G_i}), (E, h))$ gives distinct elements
in $\achow^1(X)_{\QQ}$. On the other hand,
$$\left\{
\adiff((G, \submet{h}{G}), (E, h)) \in \achow^1(X)_{\QQ} \mid
G \in \Cal{G} \right\}
$$
is finite. This is a contradiction.
So there is $G \in \Cal{G}$ with
$\adeg(\adiff((G, \submet{h}{G}), (E, h))^2) > 0$.
Thus, we get our theorem.
\qed

\subhead 3.5 Proof of Corollary~B
\endsubhead
Finally, we give the proof of Corollary~B.
We assume that
$$
\adeg\left((2 \rank E) \achern{2}{E, h} -
(\rank E - 1)\achern{1}{E, h}^2
\right) < 0.
$$
Then, by Theorem~A, there is a non-zero saturated subsheaf $F$ of $E$ with
$$
\frac{\achern{1}{F, \submet{h}{F}}}{\rank F} - \frac{\achern{1}{E, h}}{\rank E}
\in \poscone (X).
$$
Thus,
$$
\frac{\adeg(\achern{1}{F, \submet{h}{F}} \cdot A)}{\rank F} -
\frac{\adeg(\achern{1}{E, h} \cdot A)}{\rank E} > 0.
$$
This is a contradiction.
\qed


\enddocument